\documentclass[12pt]{article}
\usepackage{amsmath}
\usepackage[dvips]{graphicx}
\usepackage{amsfonts}
\usepackage{amsopn}
\usepackage{amsthm}
\title{Existence of Optimal Maps in the Reflector-type Problems}
\author{Wilfrid Gangbo\thanks{School of Mathematics, Georgia 
Institute of Technology,
Atlanta, GA 30332 (gangbo@math.\allowbreak gatech.edu).  
WG gratefully acknowledges the support of 
National Science Foundation grants DMS-00-74037, and DMS-02-00267.} 
\and Vladimir Oliker\thanks{Dept. of Mathematics and Computer Science,
Emory University, Atlanta, GA 30322, USA, 
(oliker@mathcs.emory.edu). The research of VO was partially 
supported by a grant from Emory University Research Committee and 
by the National Science Foundation grant DMS-04-05622.}
}

\date{May 4, 2005}
\def\bb{\bf} 
\newtheorem{thm}{Theorem}[section]
\newtheorem{prop}[thm]{Proposition}
\newtheorem{lemma}[thm]{Lemma}

\newtheorem{defn}[thm]{Definition}
\newtheorem{rmk}[thm]{Remark}

\def\smoothness{ {\rm (A 1)} }
\def\limite{ {\rm (A 2)} }
\def\infinite{ {\rm (A 3)} }
\def\all{ {\rm (A 1)  -- (A 3) } }

\def\proof#1{\noindent {\bf Proof{#1}:}}

\def\endproof{\ \hfill QED. \bigskip}

\def\id{\bf id}

\def\spt{{\em spt}}

\def\defd{:=}
\def\Reals{{\bb R}}
\def\Rdd{{\Reals}^{d+1} \times {\Reals}^{d+1} }
\def\Rd{{\Reals}^{d+1}}

\def\intd{\int_{\Reals^{d+1}}} 
\def\Sd{ {\bf S^d}} 

\def\calA{{\cal A}}

\def\calI{{\cal I}}

\def\calT{{\cal T}(\mu,\nu)}

\def\calP{{\cal P}}

\def\calPX{{\cal P}(\X)}
\def\calPY{{\cal P}(\Y)}

\def\spt{\text{ spt}}

\def\XY{\X \times \Y}

\def\xy{(\x,\y)}

\def\tnR{\Reals}
\def\tnRinf{\bar \Reals }

\def\X{{\bf X}}
\def\Y{{\bf Y}}

\def\a{{\bf a}}

\def\t{{\bf t}}

\def\u{{\bf u}}
\def\v{{\bf v}}
\def\x{{\bf x}}
\def\y{{\bf y}}

\def\n{{\bf n}}

\def\xy{({\bf x}, {\bf y})}

\def\0{{\bf 0}}
\def\id{{\bf id}}

\def\J{\gamma}

\def\Gam{\Gamma(\mu,\nu)}
\def\cost{{\cal C}}
\def\n{{\bf n}}

\begin{document} 
\maketitle 
\begin{abstract} 
In this paper, we
consider probability measures $\mu$ and $\nu$ on a $d$--dimensional
sphere in $\Rd, d \geq 1,$ and cost functions of the form 
$c(\x,\y)=l(\frac{|\x-\y|^2}{2})$
that generalize those arising in geometric optics where $l(t)=-\log t.$
We prove that if $\mu$ and $\nu$ vanish on $(d-1)$--rectifiable sets,
if $|l'(t)|>0,$ $\lim_{t\rightarrow 0^+}l(t)=+\infty,$ and $g(t):=t(2-t)(l'(t))^2$ is monotone then there
exists a unique optimal map $T_o$ that transports $\mu$ onto $\nu,$ where
optimality is measured against $c.$ Furthermore, $\inf_{\x}|T_o\x-\x|>0.$ Our approach is based on  direct variational arguments.
In the special case when
$l(t)=-\log t,$ existence of optimal maps on the
sphere was obtained earlier in \cite{glimm-oliker-refl2:03}
and \cite{Wang:inv03} under more restrictive assumptions. In these studies, it was assumed that  
either $\mu$ and $\nu$ are absolutely
continuous with respect to the $d$--dimensional Haussdorff measure, or they 
have disjoint supports. 
Another aspect of interest
in this work is that it is in contrast with the work in
\cite{GangboMcCann00} where it is proved that when $l(t)=t$ then
existence of an optimal map fails when $\mu$ and $\nu$ are 
supported by Jordan surfaces.
\end{abstract} 
\section{Introduction} 
\begin{figure}[ht]
\begin{center}
\includegraphics[width=12cm]{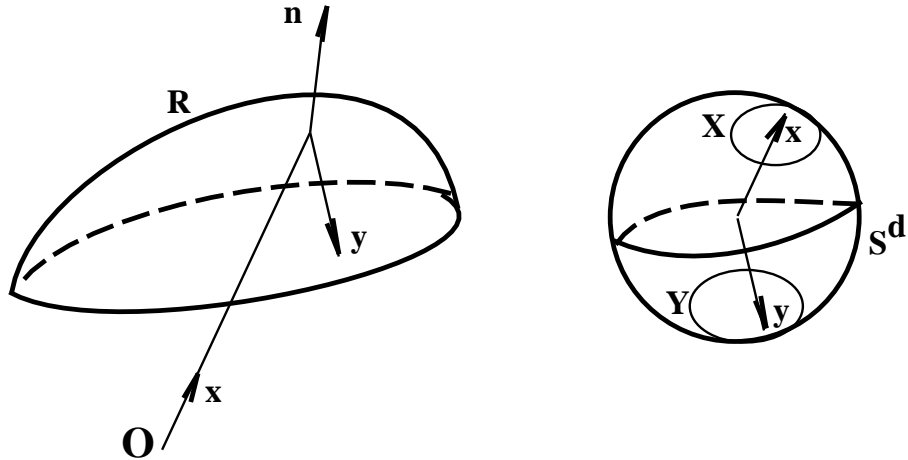} 
\caption{The reflector problem} \label{Fig. 1}
\end{center}
\end{figure}

In Euclidean space $\Rd$ consider a reflector system consisting of a point 
source $\cal{O}$ radiating with intensity $I(\x)$ in 
directions $\x \in \X$, 
where $\X$ is a closed set on a $d$--dimensional unit sphere 
$\Sd \hookrightarrow \Rd$ centered at $\cal{O}$, 
and a smooth perfectly
reflecting hypersurface $R$, star-shaped relative to $\cal{O}$, which
intercepts and reflects the light rays with directions from $\X$;
see Fig. \ref{Fig. 1}. 
Assuming  the geometric
optics approximation  and applying the classical
reflection law to determine
the set of reflected directions $\Y \subset \Sd$ (after one reflection), 
we obtain an  associated with 
$R$ ``reflector map'' $\xi: \X \rightarrow \Y$.
Assuming that no energy is lost in the process, one can apply
the energy conservation law 
to calculate the intensity distribution $L(\y)$ produced on $\Y$.
The reflector problem consists in
solving the inverse problem in which the source $\cal{O}$, 
the sets $\X, \Y$
and the intensities $I$ and $L$ are given in advance 
and the reflector $R$
needs to be determined. That is, $R$ should be such
that $\xi(\X) \supseteq \Y$ and 
$$L(\xi(\x))|J(\xi)(\x)| = I(\x)$$ for
all $\x$ in the interior of $\X$; here $J(\xi)$ denotes the Jacobian
determinant of $\xi.$

Problems of this type arise often in applications, for example,
in design of reflector antennas \cite{W}. In various forms the reflector problem has been considered by many authors and numerous papers by engineers (at least since early 60-th \cite{K} until now \cite{Wang:inv03}) are devoted to this subject. It was introduced in electrical engineering and optics independently of the mass transport problem.
Because of the strong nonlinear constraints appearing in the problem, progress has been slow and many theoretical and computational issues still remain open. The problem continues to attract considerable attention because of its practical importance and mathematical subtleties. It may be pointed out that a version of the reflector problem appears on the famous list of unsolved problems proposed
by S.T. Yau \cite{Y} in 1993.

Analytically, the reflector problem considered in this paper
can be formulated as a 
nonlinear second order elliptic partial differential equation of
Monge-Amp\`{e}re type on a subset of $\Sd$. In such form it has been studied by
V. Oliker and P. Waltman \cite{OW}, L. Caffarelli and V. Oliker \cite{CO},
X.-J. Wang \cite{Wang:inv96}, 
P.G. Guan and X.-J. Wang \cite{Pegfei_guan_wang:jdg98},
L. Caffarelli, S. Kochengin and V. Oliker \cite{CKO}, 
V. Oliker \cite{refl_geom1}, and other authors.

Recently, T. Glimm and V. Oliker \cite{glimm-oliker-refl2:03} and,
independently, X.-J. Wang \cite{Wang:inv03} have shown that 
if the function $I$ (resp. $L$) is treated as the density of a measure  
$\mu$ (resp.  $\nu$) that are absolutely continuous with respect to the 
volume measure on  $\Sd,$  then the
reflector problem can be studied as a variational problem in the
framework of Monge-Kantorovich theory, that is, a problem of
finding an optimal map
 minimizing the transport cost of transferring $\mu$ onto $\nu$ with the
cost function  $-\log(1- \x \cdot \y)$. 
In contrast with other cost functions
considered usually in the Monge-Kantorovich theory, this cost function
may assume infinite values. Consequently, in order to overcome this
difficulty, a geometric condition requiring the
supports of $\mu$ and $\nu$ to be disjoint was imposed
in  \cite{glimm-oliker-refl2:03} and \cite{Wang:inv03} to establish 
existence and uniqueness of optimal maps. Without imposing 
the condition that the supports of the measures are disjoint, existence and uniqueness of optimal maps was also obtained in \cite{glimm-oliker-refl2:03}. However, the proof is indirect as it relies on existence of weak solutions in the reflector problem established earlier in \cite{CO}, \cite{refl_geom1}.

The contribution of this study is twofold. 
First of all, we obtain existence and uniqueness 
of optimal maps $T_o$ for a class of cost functions 
that may be infinite. This class includes the logarithmic 
cost function of the reflector problem as a special case. 
The cost functions are precisely of the 
form $c(\x,\y)=l(\frac{|\x-\y|^2}{2})$ 
where $|l'(t)|>0,$ $\lim_{t\rightarrow 0^+}l(t)=+\infty,$ and $g(t):=t(2-t)(l'(t))^2$ is monotone. 
Furthermore, we prove that $\inf_{\x}|T_o\x-\x|>0,$ which we 
view as an intermediary step in the study of the regularity of 
the map $T_o.$  Secondly, our approach is variational and direct;    
the supports of measures $\mu$ and $\nu$ are allowed to overlap and it is  
merely required that these measures vanish on $(d-1)$--rectifiable 
subsets. The precise statement can be found in Theorem \ref{t1.5}.

Let us recall the main principles which ensure existence of optimal maps, which apparently, have been  explicitly pointed out for the first time in \cite{Gangbo94} and later exploited by many authors. Assume we are given a cost function $c: \Rd \times \Rd \rightarrow \Reals$ and two probability measures $\mu$ and $\nu$ on $\Rd$, say, absolutely continuous with respect to Lebesgue measure. Existence and uniqueness of an optimal map transporting $\mu$ onto $\nu$ against $c$ is ensured if  $\nabla_\x c(\x, \cdot)$ is a one-to-one map of $\Rd$ into itself. 
Note that if $c(\x,\y)=|\x-\y|^2$ then $\nabla_\x c(\x, \cdot)$ has this property. 

It is shown in \cite{GangboMcCann00} that if instead, $\mu$ and $\nu$ are 
supported by $\Sd$, then existence of an optimal map transporting $\mu$ 
onto $\nu$ against $c(\x,\y)=|\x-\y|^2$, fails. One of the points in our 
study here is that we start with a cost function $c(\x,\y)=l(|\x-\y|^2 /2)$ 
such that $\nabla_\x c(\x, \cdot)$ fails to be a one-to-one map of $\Rd$ 
into itself. However, if $\x \in \Sd$ and $\a \in \Rd$ is distinct from 
the origin, we observed that there exists a unique $\y \in \Sd$ satisfying 
the equation $\nabla_{\Sd}^\x c(\x, \y)=\a.$ Here, 
$\nabla_{\Sd}^\x c(\x, \cdot)$ is the tangential gradient of $c$ on $\Sd$
 with respect to $\x.$ The injectivity of $\nabla_{\Sd}^\x c(\x, \cdot)$ yields existence of an optimal map that transports $\mu$ onto $\nu$, if $\mu$ and $\nu$ are supported by $\Sd$ and vanish on $(d-1)$--rectifiable sets. We use the Kantorovich theory as a tool to give  a transparent explanation to the existence of a unique solution in the reflector problem, under sharp assumptions. Our main results are stated in Theorem \ref{t1.5}. We refer the reader to a variant of the reflector problem involving two reflectors
considered in a paper by T. Glimm and V. Oliker \cite{glimm-oliker-refl1:03}. We also refer the reader to a recent study by N. Ahmad \cite{Ahmad} in the plane, motivated by \cite{GangboMcCann00}. 

The paper is essentially self-contained. 
It is organized as follows. 
In order to motivate out subsequent considerations, 
we begin with a review of the reflector 
problem in section \ref{geom}. In section \ref{MK} 
we review and extend some results from the Monge-Kantorovich theory. 
Our main results establishing existence and uniqueness 
of optimal maps are in section \ref{opt}.

\section{A review of the reflector problem} \label{geom} 
Let $\X, \Y, I, L$ and $R$ be as in the introduction. 
If $\n$ is the unit normal field on 
$R$, then the incident direction $\x$ and the 
reflected direction $\y$ are related 
by the reflection law 
\begin{equation}\label{Snell}
\y=\x-2( \x \cdot\n ) \n.
\end{equation}
Thus, the hypersurface $R$ defines the {\it {reflector}} map 
$\xi:\x \rightarrow \y$ 
which maps the ``input aperture''
$\X \subset \Sd$ onto the ``output aperture''
$\Y \subset \Sd$; see Fig. 1.
The intensity
of the light reflected in direction $\y=\xi(\x)$ is given by
$
I(\x)/|J(\xi(\x))|.
$

Suppose now that the closed sets $\X$ and $\Y$ on $\Sd$ are given
as well as nonnegative integrable functions $I$ on $\X$ and $L$ on $\Y$
representing, respectively, the intensity of the source and
the desired intensity 
on the far-region $\Y$.
{\it The reflector problem 
is to determine
a reflector $R$ such that the 
map $\xi$ defined
by $R$ maps $\X$ onto $\Y$ and satisfies the equation}
\begin{equation}\label{1.2}
L(\y)=I(\xi^{-1}(\y))|J(\xi^{-1}(\y))|, \;\; \y \in \mbox{Int}\Y;
\end{equation}
see \cite{NW}, \cite{OW}. 
Note that this is an equation  on the output aperture $\Y$ 
rather
than on the input aperture $\X$. One could also set it up on 
$\X$ \cite{NO_symp}, but (\ref{1.2}) is more convenient for our
purposes here.

It was shown in \cite{NW, OW} that there exists a 
scalar quasi-potential $p: \Y \rightarrow (0, \infty)$ from which
the reflector $R$ can be recovered and in terms of which the equation 
(\ref{1.2}) when $J(\xi^{-1}) \neq 0$ is the
following equation of Monge-Amp\`ere type
\begin{equation}\label{MA}
L(\y)=I(\xi^{-1}(\y))\frac{|det[Hess(p) + (p - \rho)e]|}{\rho^n det(e)},\;
\;\; \y \in \mbox{Int}\Y,
\end{equation}
where $e$ is the 
standard metric on $\Sd$, $\rho
 = (p^2 + |\nabla p|^2)/2p$, and $Hess(p),~\nabla p$ are
 computed in the metric $e$. 
In terms of
$p$ the position vector of $R$ is given by
\begin{equation}\label{refl1}
{\bf r(\y)} = -\nabla p(\y) - (p(\y)-\rho(\y))\y: \Y \rightarrow \Rd,
\end{equation}
while
$$\xi^{-1}(\y) = \frac{{\bf r}(\y)}{\rho(\y)}.$$
Note that $|{\bf r}| = \rho$.

A close examination of  (\ref{refl1}) shows that
it describes $R$ as an envelope of a family of 
paraboloids of revolution $P(\y)$ tangent to $R$,
parametrized by their axes
$\y \in \Y$ and with common focus $\cal{O}$. 
For each $\y,$ $p(\y)$ is the focal parameter of $P(\y)$.
This observation was used in \cite{CO} for
the weak formulation of the reflector problem where a class
of convex reflectors corresponding to positively elliptic solutions
of (\ref{MA}) was introduced. Reflectors corresponding to negatively
elliptic solutions of (\ref{MA}) can be introduced and analyzed
by similar methods \cite{glimm-oliker-refl2:03}, \cite{Wang:inv03}. 
Such reflectors, however, are only piecewise concave
(relative to the origin $\cal{O}$). For brevity we discuss here only
convex reflectors which we now define.

 Let $\Y$ be a subset on $\Sd$ and $p^{\prime}:\Y \rightarrow (0,\infty)$
a bounded function. Let   $\{P(\y)\}$ be a family 
of paraboloids of revolution
with axes of direction $\y \in \Y$,
common focus $\cal{O}$ and polar radii
\begin{equation}\label{polar1}
\rho_\y(\x)=\frac{p^{\prime}(\y)}{1-\x \cdot \y }, \;\; 
\x \in \Sd \setminus \{ \y\}.
\end{equation} 
The closed convex hypersurface $R$ given by 
${\bf r}(\x) = \rho(\x)\x, ~\x \in \Sd,$ with 
\begin{equation}\label{polar2}
\rho(\x) = \inf_{\y\in \Y}\rho_y(\x),~\x \in \Sd,
\end{equation} 
is called a {\bf {reflector defined by the family
$\{P(\y)\}_{\y \in \Y}$}} (with the light source at $\cal{O}$).

Let $R$ be a reflector, $z$ a point on $R$,
$P$ a paraboloid of revolution with focus at $\cal{O}$ 
and $B$ the convex
body bounded by $P$. If $R \subset B$ and
$z \in P \bigcap R$ then   $P$ is called
{\bf supporting to $R$
at $z$}.

For any $\y \in \Sd$ a reflector $R$ has a supporting paraboloid 
with axis $\y$.
The corresponding continuous function giving
the focal parameters of all such 
supporting paraboloids $P(\y),~\y \in \Sd$, of $R$ is called the
{\it{focal function}} and denoted by $p$. The natural question
of characterization of focal functions of closed convex reflectors was 
partially answered in \cite{refl_geom1}.

The {\bf {reflector  map}} (possibly multivalued)
 generalizing (\ref{Snell}) and denoted again
by $\xi$ is defined for $\x \in \Sd$ as
\begin{equation}\label{raytr}
\xi(\x)=\{ \y\in \Sd\;\bigl|\; \text{$P(\y)$
is supporting to $R$ at $\rho(\x) \x$}\}.
\end{equation} 
It follows from (\ref{polar1}) and (\ref{polar2}) that, 
equivalently, the reflector map can be defined as
\begin{equation}\label{raytr1}
\xi(\x)=\{ \y\in \Sd\;\bigl|\; \log\rho(\x) - \log p(\y) 
= -\log (1- \x\cdot \y)\}.
\end{equation}
The total amount of energy transferred from the source $\cal{O}$
 in a given
set of directions $\omega \subset \Y$ is best described with the use
of   the map $\xi^{-1}$ defined on any subset $\omega \subset \Y$ by setting
\[
\xi^{-1} (\omega) = \bigcup_{\y\in\omega}\xi^{-1}(\y).
\]
It is shown in \cite{CO, refl_geom1} that for any Borel subset
$\omega \subset \Y$ the set $\xi^{-1}(\omega)$ is Lebesgue measurable 
on $\Sd$.
Thus, if $I \geq 0$ is the intensity of the source then
the total amount of energy transferred by $R$
to the set $\omega \subset \Y$  is given by the ``energy'' function
\begin{equation} \label{energy}
G(R,\omega )=\int_{\xi^{-1}(\omega)} I(\x) d \sigma(\x),
\end{equation}
where $d\sigma$ is the standard $d-$volume form on $\Sd$. It is 
assumed here
and everywhere below that $I$ is extended from $\X$ 
to the entire $\Sd$
by setting $I(\x)\equiv 0~~\forall \x \in \Sd \setminus \X$. 
The function $G(R,\omega )$ is a Borel measure on $\Y$
(not necessarily absolutely continuous)  \cite{CO}, \cite{refl_geom1}.

For a given nonnegative and integrable function $L$ on $\Y$ 
we say that a closed convex reflector $R$ is a {\bf weak solution}
 of the reflector
problem if
\begin{equation}\label{reqn}
G(R,\omega ) = \int_{\omega}L(\y)d\sigma(\y) ~~~\mbox{for any Borel set}
~~~ 
\omega \subset \Y.
\end{equation}

Put 
\begin{equation}\label{def-mu-nu}
\mu[A]= \int_A I(\x) d\sigma(\x), \qquad \nu[B]= \int_B L(\y) d\sigma(\y)
\end{equation}
for $A \subset \X$ and $B \subset \Y$ Borel sets.
An obvious necessary condition for existence of a weak solution to the
reflector problem is that the total energy of the source and the
total energy on the output aperture are in balance, that is,
\begin{equation}\label{bal0}
\mu[\X]= \nu[\Y].
\end{equation}

Excluding the trivial case when either $\mu[\X]$ or $\nu[\Y]$ is
zero, it may be assumed, without loss of generality, 
that measures satisfying (\ref{bal0}) are normalized so that
they are probability measures.

The following existence result was established by L. Caffarelli and V. Oliker
in \cite{CO} (and reproduced  partly in \cite{CKO} and \cite{refl_geom1}).
\begin{thm}\label{reflector_existence} 
Let $\X$ and $\Y$ be closed sets on $\Sd$
(possibly coinciding with $\Sd$) and
$I\geq 0 $, $L\geq 0$  two probability densities on $\X$ and $\Y$,
respectively.
Then there exists a reflector which is a weak solution of the
reflector problem.
\end{thm}
The proof 
is obtained in two steps. First the problem
is solved in the case when the right hand side
in (\ref{reqn}) is a finite
sum  of Dirac masses. In this case a constructive
minimization procedure together with an apriori two-sided $C^0$ estimate
of $\rho$ is used to obtain the weak solution, which
is also unique. The general problem is solved
by approximating the right hand side in
 (\ref{reqn}) by finite sums of Dirac masses and obtaining the
solution $R$ as a limit of a sequence of special solutions
$R_k, k =1,2,...,$ constructed on
the previous step. The measures $G(R_k, \cdot)$ are weakly
continuous and converge to $G(R, \cdot)$. Consequently, $R$ is indeed
a weak solution of the reflector problem. 

The procedure we have just described proves, in fact,
 that the reflector problem 
admits a solution if the right hand side
in (\ref{reqn}) is any nonnegative probability measure on $\Y$, 
possibly with a singular part.

Existence of regular solutions was studied by X.-J. Wang \cite{Wang:inv96}
and P.-F. Guan and X.-J. Wang \cite{Pegfei_guan_wang:jdg98}. 

In the
framework of Monge-Kantorovich theory the reflector problem was studied
by T. Glimm and V. Oliker \cite{glimm-oliker-refl2:03} and
X.-J. Wang \cite{Wang:inv03}. The following result was proved in 
\cite{glimm-oliker-refl2:03}, Theorem 4.1.
\begin{thm} \label{refl_minimizer} Let
$R$ be a weak solution of the reflector problem with the
reflector map $\xi$. 
 Then  $\xi_\# \mu=\nu$ (that is, $\xi$ pushes $\mu$ forward to $\nu$)
and it is a minimizer of the problem 
\begin{equation}\label{e05.11.1}
\inf_{\xi'} \{\int_\X -\log(1- \x \cdot \xi'(\x))d \mu(\x) 
\; | \; \xi'_\# \mu=\nu \}.
\end{equation}
Furthermore, any other minimizer of (\ref{e05.11.1}) 
is equal to $\xi$ almost everywhere 
on the set $\{\x \in \Sd ~\bigl |~I(\x) \not = 0\}.$
\end{thm} 
\begin{rmk} 
Theorems \ref{reflector_existence} and \ref{refl_minimizer}, together,
imply
 existence of minimizers to the problem
(\ref{e05.11.1}). In addition, if $I > 0$ on $\X$ and
$\X$ is connected, 
these results imply 
that, except for a set of measure zero, 
any minimizer of (\ref{e05.11.1})  (with $\mu$ and $\nu$ as in
(\ref{def-mu-nu})) is
a reflector map associated 
with a closed convex reflector
in $\Rd$.  Such reflector is unique up to a  constant multiple
 of the function $\rho(\x)$
in (\ref{polar2}); see \cite{glimm-oliker-refl2:03}. 
\end{rmk}

On the other hand, the minimization problem
(\ref{e05.11.1}) is a variant of the Monge problem on $\Sd$ 
(see the beginning of section \ref{MK}, below).
By a different method, in the framework of Monge-Kantorovich theory, 
the existence and uniqueness of minimizers to (\ref{e05.11.1}) 
was proved in \cite{glimm-oliker-refl2:03} and \cite{Wang:inv03} 
under the additional condition 
that $\spt (\mu) \cap \spt (\nu) = \emptyset$. 
Furthermore, if  $\rho$ and $p$ are the functions in (\ref{raytr1}) 
then $u_o(\x)= \log \rho(\x)$ and $v_o(\y)= -\log p(\y)$ maximize 
$$
(u,v) \rightarrow \int_\X u(\x) d \mu(\x) +\int_\Y v(\y) d \nu(\y)
$$
over the set of pairs $(u,v) \in C(\X \times \Y)$ 
satisfying $u(\x) +v(\y) \leq  -\log(1- \x \cdot \y)$ for 
all $(\x,\y) \in \X \times \Y.$  
\section{Background on the Monge-Kantorovich theory}\label{MK}
As observed in section \ref{geom},  
the reflector problem can be stated as a variant of the problem 
of Monge with the 
cost function $c \xy = -\log(1-\x \cdot \y),~ \x, \y \in \Sd$. 
In this section we recall some facts of the 
Monge-Kantorovich theory that will be needed  to study 
an analogue of the reflector 
problem with general cost functions.

We first fix some notation. For a set $Z \subset \Rd$ we denote
by $\calP(Z)$ the
set of Borel probability measures on $Z.$ As usual,
if $Z$ is a closed subset of $\Rd$ and $\gamma \in
\calP(Z)$ then the support of $\gamma$ is the smallest 
closed set $\text{ spt} \gamma
\subset Z$ such that $\gamma[\text{ spt} \gamma]= \gamma[Z]=1.$

Suppose, $\X, \Y \subset \Rd$ are closed sets. If $\mu
\in \calP(\X),$ $\nu \in \calP(\Y),$ we denote by $\Gam$ the set of
joint measures $\gamma$ on $\Rdd$ that have $\mu$ and $\nu$ as their
{\em marginals}:  $\mu[U] = \gamma[U \times \Rd]$ and $\gamma[\Rd
\times U]= \nu[U]$ for Borel $U \subset \Rd.$ In fact, if $\gamma \in
\Gam$ then $\spt \gamma \subset \X \times \Y.$

Assume that we are given two probability measures
$\mu$ and $\nu$ on $\Rd.$ Let $\Gam$ be the set of  joint measures
$\gamma$ on $\Rdd$ that have $\mu$ and $\nu$ as their marginals.
Kantorovich's problem is to minimize the {\em
transport cost} 
\begin{equation}\label{cost} 
{\cost(\gamma) \defd \int
c \xy \, d\gamma \xy } 
\end{equation} 
for some given $c$ among joint measures $\gamma$  in
$\Gam$, to obtain 
\begin{equation}\label{Kantorovich} {\inf_{\J \in
\Gam} \cost[\J]}.  
\end{equation} 

Let $\calT$ be the set of Borel maps
$T: \Rd \rightarrow \Rd$ that push $\mu$ forward to $\nu:$
$\mu[T^{-1}(B)]=\nu[B]$ for all Borel sets $B \subset \Rd$.  The Monge
problem is to minimize $$ \calI[T] = \intd c(\x, T\x) d \mu(\x) $$ over
the set $\calT.$ 

There is a natural embedding which associates to $T \in \calT $ a $\gamma_T:= (\id  \times T)_\# \mu \in \Gam,$ where $\id: \Rd \rightarrow \Rd$ is the identity map. Since $\cost[\gamma_T]=\calI[T]$ we conclude that
\begin{equation}\label{relaxation} {\inf_{\J \in \Gam} \cost[\J]} \leq
{\inf_{T \in \calT} \calI[T]}.  
\end{equation}

Throughout this section we use the notation 
$\Reals \cup \{+\infty\}=\tnRinf$ and assume that 
$c: \Rdd \rightarrow \tnRinf.$ We endow $\tnRinf$ with the usual 
topology so that $c \in C(\Rdd,\tnRinf)$ means that  
$$
\lim_{(\x, \y) \rightarrow (\bar \x,\bar \y)} c( \x, \y)=c(\bar \x,\bar \y).
$$ 
In particular, if $c(\bar \x,\bar \y)=+\infty$ then $c( \x, \y)$ 
tends to $+\infty$ as $(\x, \y)$ 
tends to $(\bar \x,\bar \y).$
\begin{defn}\label{d1.1} A subset $S \subset \Rdd$ is said to be 
$c$--cyclically 
monotone if for every natural number $n$ and every $\{(\x_i, \y_i) \}_{i=1}^n 
\subset S$ we have  
$$
\sum_{i=1}^n c(\x_i , \y_i) \leq \sum_{i=1}^n c(\x_i , \y_{\sigma(i)} ),
$$
for all permutation $\sigma$ of $n$ letters. 
\end{defn}
This notion of $c$--cyclical monotonicity was introduced 
by Knott and Smith 
\cite{KnottSmith84}. 
When $c(\x,\y)=|\x-\y|^2$, $c$--cyclical monotonicity 
is simply called 
cyclical monotonicity \cite{Rockafellar}. 
%

\begin{prop}\label{p1.2} Assume that $\X, \Y \subset \Rd$ are closed sets, 
that 
$\mu \in \calPX,$  $\nu \in \calPY$ and that $c \geq 0$ is lower 
semicontinuous on $\XY.$ Then, 

(i) there is at least one optimal measure $\gamma_o \in \Gam$. 

(ii) Suppose that in addition $c \in C(\XY, \tnRinf).$  
Unless $\cost \equiv +\infty$ throughout $\Gam,$ 
there is a $c$--cyclically monotone set $S \subset \Rdd$ containing 
the support of all optimal measures in $\Gam.$
\end{prop}
%

Part (i) of proposition \ref{p1.2} can be found in \cite{Kellerer84} 
(Theorem 2.19). Note that that proof is interesting only in the case $\cost \not \equiv +\infty.$ Part (ii) was established in \cite{AbdellaouiHeinich94} 
and   \cite{GangboMcCann96} for $c \in C(\XY, \tnR).$  One can readily 
adapt the proof of Theorem 2.3 and corollary 2.4 of \cite{GangboMcCann96} 
to cost functions  $c \in C(\XY, \tnRinf).$ 

%

\begin{defn}\label{d1.3} Suppose that 
$\psi: \Rd \rightarrow \Reals \cup \{-\infty\}$ is not identically $-\infty.$ 
Then (i) $\psi$ is said to be $c$--concave if there exists a set $\calA 
\subset \Rd 
\times \Reals$ such that 
\begin{equation}\label{e1.66}
\psi(\x)= \inf_{(\y, \lambda) \in \calA} c \xy + \lambda, \quad \quad (\x \in 
\X).
\end{equation} 
(ii) The $c$--superdifferential $\partial^c \psi$ of $\psi$ consists of the 
pairs $ \xy \in \Rdd$  for which  
\begin{equation}\label{e1.6}
c \xy - \psi(\x) \leq c(\v, \y) -\psi(\v),
\end{equation}
for all $\v \in \Rd.$ 
\hfill\break 
(iii) We define $\partial^c \psi(\x) \subset \Rd$ to be the set of $\y$ such 
that $(\x, \y) \in \partial^c \psi(\x).$ If $E \subset \Rd$, $\partial^c 
\psi(E)$ is the union of the $\partial^c \psi(\x)$ such that $\x \in E.$
\hfill\break
(iv) The c--transform of $\psi$ is the function 
$\y \rightarrow \psi^c(\y) = \inf_{\x \in \X}\{c(\x,\y) -\psi(\x)\}.$
\end{defn}
\begin{rmk}\label{r1.3bis} Suppose that 
$\psi: \Rd \rightarrow \Reals \cup \{-\infty\}$ is not identically $-\infty$ 
and is given by (\ref{e1.66}).  We have

(i) $\psi^c(\y) \geq -\lambda>-\infty$ if $(\y, \lambda) \in \calA$ where $\calA$ is the set in (\ref{e1.66}). Hence, $\psi^c \not \equiv -\infty.$

(ii) $\psi^{cc}=\psi.$
\end{rmk}
\proof{} The proofs of these remarks are well documented 
when $c: \Rdd \rightarrow \Reals.$ 
We verify that the same proofs apply when $c$ may take the value $+\infty.$ 

If  $(\y, \lambda) \in \calA$ then $-\lambda \leq c(\x,\y) -\psi(\x)$ 
for all $\x \in \X.$ Hence $\psi^c(\y)$, the infimum of $c(\x,\y) -\psi(\x)$ 
over $\X$, is not smaller than $-\lambda.$ This proves (i). 

The inequality $\psi^{cc} \geq \psi$ which holds for general functions is 
readily checked. It remains to prove that when (\ref{e1.66}) holds, 
then $\psi^{cc} \leq \psi.$ Fix $\x \in \X$ and 
let $\{ (\y_n, \lambda_n)\}_{n=1}^\infty \subset \calA$ be such that 
$$
\psi(\x)= \lim_{n \rightarrow +\infty} c(\x, \y_n) + \lambda_n.
$$
By (i), $\lambda_n \geq -\psi^c(\y_n)$ and so, 
$$
\psi(\x) \geq \limsup_{n \rightarrow +\infty} c(\x, \y_n)-\psi^c(\y_n) \geq 
\psi^{cc}(\x).
$$
This proves (ii). 
\endproof

\section{Existence and uniqueness of optimal maps} \label{opt}
Throughout this
section, we assume that 
\smallskip \\ {\em \smoothness\ $c \in C^1(\Rdd
\setminus \Delta).$} 
 \smallskip \\ {\em \limite\ $c:\Rdd \rightarrow [0,+\infty]$ is lower
semicontinuous.  } 
 \smallskip \\ {\em \infinite\ for any $\x_o \in \Rd$, $c(\x_o,
\x_o)=+\infty.$  } 
\hfill\break
\hfill\break
We are interested in probability measures $\mu,$  $\nu$ for which 
$\cost \not \equiv +\infty$ throughout $\Gam.$ Assume that 
$c \in C(\Rdd \setminus \Delta, \Reals),$ where 
$\Delta:=\{(\x,\x) \; | \; \x \in \Rd \}$ 
denotes the diagonal. Proposition \ref{p03.08.2} 
provides a sufficient condition which ensures that 
$\cost \not \equiv +\infty$ throughout $\Gam.$ 
Before stating that proposition, let us introduce the sets 
$$
S(a,b)=\{\x=(x_1, \cdots, x_d,x_{d+1}) \; | \;  a \leq x_{d+1} \leq b \}, 
$$
for $-1 \leq a \leq b \leq 1.$  

\begin{prop} \label{p03.08.2} Suppose that 
$\mu, \nu \in \calP(\Sd)$ are Borel measures which 
vanish on $(d-1)$--rectifiable sets. Then 
there exists $\gamma \in \Gam$ and $\epsilon>0$ such that 
$$
|\x-\y| \geq \epsilon 
$$
for all $(\x,\y) \in \spt \gamma.$
\end{prop}
\proof{} Since $\mu$ and $\nu$ vanish on $(d-1)$--rectifiable sets, 
the functions 
\begin{equation}\label{e03.09.1}
t \rightarrow \mu[S(a,t)], \quad t \rightarrow \nu[S(a,t)] 
\quad \hbox{are continuous}.
\end{equation} 

{\bf Case 1.} Assume first that there exists $c \in (-1,1)$ such that 
\begin{equation}\label{e03.08.1}
\spt \mu \subset S(-1,c), \quad \spt \nu \subset S(c,1).
\end{equation}
Thanks to (\ref{e03.09.1}) we may choose $\epsilon_1, \epsilon_2 >0$ such that 
$$
\mu[S(-1,-1 + \epsilon_1)]= \nu[S(c,c + \epsilon_2)]= 1/2.
$$ 
By (\ref{e03.08.1}) $-1+\epsilon_1<c.$ Set
$$
\bar \gamma= 2( \mu^- \otimes \nu^- + \mu^+ \otimes \nu^+)
$$
where 
$$
\mu^-= \mu|_{S(-1,-1+\epsilon_1)}, \quad \mu^+= \mu|_{S(-1+\epsilon_1,c)}, 
\quad 
\nu^-= \nu|_{S(c,c+\epsilon_2)}, \quad \nu^+= \nu|_{S(c+\epsilon_2,1)} . 
$$
Note that $\bar \gamma \in \Gam$ and 
$$
|\x-\y| \geq \min\{ c+1-\epsilon_1, \epsilon_2\}>0 
$$
for all $(\x,\y) \in \spt \gamma.$ This proves the proposition in this 
special case.
\hfill\break

{\bf Case 2.} Assume that $\spt \mu$ and $\spt \nu$ are arbitrary. 
We use (\ref{e03.09.1}) to choose $c \in (-1,1)$ such that 
\begin{equation}\label{e03.08.2}
\mu [ S(c,1)]= \nu [ S(-1,c)]:=m.
\end{equation}
If $m=0$ then 
$$
\spt \mu \subset S(-1,c), \quad \spt \nu \subset S(c,1)
$$
and so, we reduce the discussion to the case 1. Similarly, if 
$m=1$ we reduce the discussion to the case 1. 

Assume in the sequel that $0<m<1.$ Set 
$$
\mu^-=\mu|_{S(-1,c)}, \quad \mu^+=\mu|_{S(c,1)}, 
\quad \nu^-=\nu|_{S(-1,c)}, \quad \nu^+=\nu|_{S(c,1)}.
$$
By (\ref{e03.08.2}), $\frac{\mu^+} {m} $ and $\frac{\nu^-}{m} $ 
are probability measures. They satisfy 
$$
\spt \nu^- \subset S(-1,c), \quad \spt \mu^+ \subset S(c,1).
$$
Having that $\nu^-,$ $\mu^+$ satisfy the assumptions of case 1, 
we may find $\bar \gamma \in \Gamma(\frac{\mu^+} {m} , \frac{\nu^-} {m})$ 
and $\bar \epsilon >0$ such that 
\begin{equation}\label{e03.08.3}
|\x-\y| \geq \bar \epsilon
\end{equation}
for all $(\x,\y) \in \spt \bar \gamma.$ Similarly, there exists 
$\tilde \gamma \in \Gamma(\frac{\mu^-} {1-m} , \frac{\nu^+} {1-m})$ 
and $\tilde \epsilon >0$ such that 
\begin{equation}\label{e03.08.4}
|\x-\y| \geq \tilde \epsilon
\end{equation}
for all $(\x,\y) \in \spt \tilde \gamma.$ 

Set 
$$
\gamma= m \bar \gamma + (1-m) \tilde \gamma.
$$
Then $\gamma \in \Gam$ and by (\ref{e03.08.3})--(\ref{e03.08.4}), we have that 
$$
|\x-\y| \geq \min\{\bar \epsilon, \tilde \epsilon \}>0
$$
for all $(\x,\y) \in \spt \tilde \gamma.$ This proves the proposition. 
\endproof
%
%
%
%
%
\begin{lemma} \label{p1.2ter} Suppose that $c$
satisfies \all and that $\X, \Y \subset \Rd$ are closed sets. Suppose
that $S \subset \XY$ is $c$-cyclically monotone and contains two pairs
$(\x_o, \y_o)$, $(\bar \x_o, \bar \y_o)$ such that $\x_o\not =\bar \x_o$ 
and $\y_o \not= \bar \y_o .$ Then, there exists a function $F \in
C(\XY)$ depending only on the pairs such that
\begin{equation}\label{e1.2.1} c( \x_o, \y_o) + c(\bar \x_o, \bar \y_o)
+ c(\x,\y) \leq F(\x,\y) 
\end{equation} 
for all $(\x,\y) \in S.$
\end{lemma} 
 \proof{} Define 
\begin{equation}\label{e1.2.2}
\begin{cases} F(\x,\y) &= \min\{ c(\x, \y_o) + R_1(\y), c(\x, \bar
\y_o) + R_2(\y) \} \cr R_1(\y) &= \min\{ c(\x_o, \y) + c(\bar \x_o,
\bar \y_o), c(\x_o, \bar \y_o) + c(\bar \x_o, \y)\} \cr R_2(\y) &=
\min\{ c(\x_o, \y) + c(\bar \x_o, \y_o), c(\x_o, \y_o) + c(\bar \x_o,
\y) \} \cr 
\end{cases} 
\end{equation} 
We use that $\x_o \not =\bar
\x_o$ and \all to obtain that $R_1, R_2 \in C(\Y).$ This, together with
the fact that $\y_o \not =\bar \y_o$ yields that $F$ is continuous on
$\XY.$ If $(\x,\y)$ is another element of $S$ then setting $$ (\x_1,
\y_1)=(\x_o, \y_o), \quad (\x_2, \y_2)=(\bar \x_o, \bar \y_o),  \quad
(\x_3, \y_3)=(\x, \y) $$ and using the $c$--cyclical monotonicity of
$S$, we obtain (\ref{e1.2.1}).  
\endproof 

It is well known that a set is cyclically monotone if and only if it is
contained in the subdifferential of a convex function
\cite{Rockafellar}.  An analogue of this result was proved by Smith and
Knott \cite{SmithKnott92} for general cost functions $c: \XY
\rightarrow \Reals$.  The following Lemma \ref{p1.4} is a further
extension  that is needed to deal with cost functions satisfying \all
(and may be $=+\infty$ somewhere).  Below, we check first that the
proof in \cite{Ruschendorf96} extends to such cost functions  and then
we show that the infimum in (\ref{e1.66ter}) can be performed over the
subset of $\y$ such that $|\x -\y|> \delta > 0.$ 
\begin{lemma}\label{p1.4} Suppose that $c$ satisfies \all and  $\X, \Y
\subset \Rd$ are compact sets. Suppose that $S \subset \XY$ is
$c$-cyclically monotone and contains two pairs $(\x_o, \y_o)$, $(\bar
\x_o, \bar \y_o)$ such that $\bar \x_o, \bar \y_o \not \in \{\x_o, \y_o\}$ Then, \hfill\break (i) $S$ is contained in the
$c$--superdifferential of a $c$--concave function $\psi: \Rd
\rightarrow \Reals \cup \{-\infty\}$ such that there exists $\delta>0$
satisfying 
\begin{equation}\label{e1.66ter} \psi(\x)= \inf_{\y }
\Bigl\{ c(\x,\y) - \psi^c(\y) \; | \; |\x -\y| \geq \delta \Bigr\},
\quad \quad (\x \in \X).  
\end{equation} 
\hfill\break (ii) If $(\x, \y)
\in \partial^c \psi$ then $\psi(\x)>-\infty$ and $|\x-\y| \geq \delta$
for some $\delta>0$ that depends only on $\psi.$ 
\end{lemma} 
\proof{} The expression $\psi(\x)= \inf_{\y \in \Y} \Bigl\{ c(\x,\y) -
\psi^c(\y)\Bigr\},~\x \in \X,$ is well-known in the literature. The
only new and useful fact we want to point out is that $\psi(\x)$ will
be obtained by minimizing $c(\x,\y) - \psi^c(\y)$ not on $\Y$ (as it is
usually done), but on $\{y \in \Y \; | \; |\x -\y| \geq \delta \}.$ For
completeness, we give the detailed proof below.

Since $\XY$ is compact, the function $F$ defined in (\ref{e1.2.1})
attains its maximum. We use Lemma \ref{p1.2ter} and the fact that $c$
is lower semicontinuous and equals $+\infty$ on $\Delta$ to conclude
the following: if $S' \subset \XY$ is $c$--cyclically monotone and
contains $(\x_o, \y_o)$, $(\bar \x_o, \bar \y_o)$ then there exists
$\delta>0$ such that $|\x -\y| \geq 2 \delta$ for every $(\x, \y) \in
S'.$ In particular, $\x \not =\y$ for $(\x, \y) \in S'.$ In particular, there exists $\delta_S>0$ such that $|\x -\y| \geq 2 \delta_S$ for every $(\x, \y) \in
S.$

As in \cite{Ruschendorf96}, we define 
\begin{equation}\label{e1.66bis}
\psi(\x)= \inf_{n} \inf_{\{(\x_i,\y_i)\}_{i=1}^n \subset S}
\biggl\{c(\x,\y_n) + \sum_{j=0}^{n-1}  c(\x_{j+1}, \y_{j}) -
\sum_{j=0}^{n} c(\x_j, \y_{j})\biggr\}.  
\end{equation} 
Since $c$ is
finite on $S$ and nonnegative on $\Rdd$, we conclude that $\psi(\x)$ is
well-defined.

The $c$--cyclical monotonicity of $S$ gives that $$ c(\x_o, \y_n)+
\sum_{j=0}^{n-1}  c(\x_{j+1}, \y_{j}) - \sum_{j=0}^{n}c(\x_j, \y_{j})
\geq 0 $$ and so, $\psi(\x_o) \geq 0.$ Taking $n=1,$ $\x_1=\x_o$ and
$\y_1=\y_o$ in (\ref{e1.66bis}) gives that $\psi(\x_o) \leq 0.$ We
conclude that $\psi(\x_o) = 0$ and so, $\psi$ is not identically
$-\infty.$ This, together with the fact that  $\psi$ is clearly of the
form (\ref{e1.66}) yields that $\psi$ is $c$--concave.

\hfill\break {\bf Claim 1.} We have that $S \subset \partial^c \psi.$

Fix $(\x, \y) \in S$ and for each $m$ integer, let
$\{(\x_i,\y_i)\}_{i=1}^{n_m} \subset S$ be such that $$ \lim_{m
\rightarrow +\infty} \psi_m =\psi(\x), $$ where $$
\psi_m:=c(\x,\y_{n_m}) + \sum_{j=0}^{{n_m}-1}  c(\x_{j+1}, \y_{j}) -
\sum_{j=0}^{{n_m}} c(\x_j, \y_{j}).  $$ Setting $(\x_{{n_m}+1}, 
\y_{{n_m}+1})=(\x,\y)$ we have that 
$$
-c(\x,\y)+\psi_m=c(\v, \y_{{n_m}+1}) -c(\v, \y) + \sum_{j=0}^{{n_m}}
c(\x_{j+1}, \y_{j}) - \sum_{j=0}^{{n_m}+1} c(\x_j, \y_{j}) \geq
\psi(\v)-c(\v, \y).  
$$ 
Letting $m$ go to $+\infty$ we conclude that
\begin{equation}\label{efeb17.1} 
\psi(\x)-c(\x, \y) \geq \psi(\v)-c(\v,\y),  
\end{equation} 
which proves claim 1.

\hfill\break {\bf Claim 2.} Whenever $(\x,\y) \in \partial^c \psi,$ we 
have that $\psi(\x)>-\infty$  and $\x \not = \y.$ 
 
 Recall that $(\x,\y) \in \partial^c \psi$ is equivalent to (\ref{efeb17.1}). 
Setting $(\x,\y)=(\bar \x_o, \bar \y_o),$ $\v=\x_o$ in (\ref{efeb17.1}), 
using the facts that $\psi(\x_o)=0$ and $\x_o \not = \bar \y_o$ we 
obtain that $\psi(\bar \x_o)$ is finite. Next, if $ (\u_o ,\v_o) \in S,$ 
setting $\v=\u_o$ we have that 
\begin{equation}\label{e05.11.4} 
c(\x, \y) -\psi(\x) \leq c(\u_o, \y) -\psi(\u_o),   
\end{equation}
If $\y \not = \x_o$ we set $\u_o= \x_o$ in (\ref{e05.11.4}) to obtain the claim. If $\y = \x_o$, we set $\u_o= \bar \x_o$ and we use the fact that $\psi(\bar \x_o)$ is finite to obtain the claim. 

\hfill\break {\bf Claim 3.} The set $\partial^c \psi$ is $c$--cyclically 
monotone. 

For the sake of completeness, we reprove this claim although it is a 
repetition of a known argument in \cite{Ruschendorf96}. 
If $\{(\x_i,\y_i) \}_{i=1}^n \subset \partial^c \psi,$ 
setting $(\x_{n+1},\y_{n+1})=(\x_1,\y_1)$ we have that  
\begin{equation}\label{e05.11.7} 
c(\x_i,\y_i) -\psi(\x_i) \leq c(\x_{i+1},\y_i) -\psi(\x_{i+1}).  
\end{equation}
By claim 2, each term in (\ref{e05.11.7}) is finite and so, 
$$  
0= \sum_{i=1}^n \psi(\x_{i+1})-\psi(\x_i) \leq \sum_{i=1}^n c(\x_{i+1},\y_i) -c(\x_i,\y_i),  
$$ 
which proves the claim.
\hfill\break 

We use Lemma \ref{p1.2ter}, the facts that 
$\partial^c \psi$ is $c$--cyclically monotone, 
that $(\x_o, \y_o),$ $(\bar \x_o, \bar \y_o) \in
\partial^c \psi,$ that $\x_o \not= \bar \x_o$ and 
$\y_o \not= \bar \y_o,$ to obtain the existence of 
some $\delta>0$ such that $|\x-\y| \geq 2 \delta$ 
for all $(\x,\y) \in \partial^c \psi.$ By remark \ref{r1.3bis} 
(ii) $\psi = (\psi^c)^c$ and so if $\x \in \X,$ there exists 
a sequence $\{\y_n \}_{n=1}^\infty \subset \Y$ such that 
\begin{equation}\label{efeb17.4} \psi(\x) = \lim_{n
\rightarrow +\infty} c(\x,\y_n) - \psi^c(\y_n).  \end{equation} Since
$\Y$ is compact, we may extract from $\{\y_n \}_{n=1}^\infty$ a
subsequence (which we still label $\{\y_n \}_{n=1}^\infty$) that
converges to some $\y \in \Y.$ Recall that  the function $\psi^c$ is
upper semicontinuous as an infimum of upper semicontinuous functions
and therefore, (\ref{efeb17.4}) yields 
$$
\psi(\x) \geq c(\x, \y) -\psi^c(\y).  
$$ 
This proves that  $\psi(\x) = c(\x, \y) -\psi^c(\y)$ and so,  
$(\x,\y) \in \partial^c \psi$. Hence $|\x-\y| \geq 2 \delta$ 
and $\y$ is a minimizer in (\ref{e1.66ter}). \endproof 

\begin{rmk} \label{r1.2bis} Suppose that $\mu \in \calPX$ and $\nu \in
\calPY$ have no atoms and that $\gamma$ minimizes $\cost$ over $\Gam$
and that $\cost(\gamma)<+\infty.$ Then, $\gamma[\Delta]=0$ and so,
$\spt(\gamma) \setminus \Delta$ contains at least one element, say
$(\x_o, \y_o).$ Also, $\gamma[E]=0$ where $$ E= \Bigl(\{ \x_o, \y_o \}
\times \Y \Bigr)  \cup \Bigl( \X \times \{ \x_o, \y_o \} \Bigr) $$
Hence, the set $\XY \setminus (E \cup \Delta))$ is nonempty, and so, it
contains an element $(\bar \x_o, \bar \y_o).$ Note that 
$\bar \x_o, \bar \y_o \not \in \{\x_o, \y_o \}.$  
\end{rmk} 

We now further specialize the set of cost functions under
consideration by assuming that \begin{equation}\label{e1.5.1} c(\x,\y)=
\begin{cases} l(\frac{ |\x-\y|^2} {2}) & \quad \text{if} \quad \x \not
=\y \cr +\infty & \quad \text{if} \quad \x =\y \cr \end{cases}
\end{equation} where $l \in C^2(0,+\infty)$ is such that
\begin{equation}\label{e1.5.2} \begin{cases} \lim_{t \rightarrow 0^+}
l(t) & = +\infty \cr |l'(t)| & >0 \quad \quad (t>0) \end{cases}
\end{equation} 
Define 
\begin{equation}\label{e03.10.1}
g(t)=t(2-t)(l'(t))^2 
\end{equation} 
Note that if $g$ is monotone on
$(0,2]$ then $g[\frac{\delta}{2},2]$ is a closed interval on which
$g^{-1}$ exists and is continuous. Define 
$$
M(\a,\x)=\Bigl[1-g^{-1}(|\a|^2) \Bigr] \x -
 \frac{\a}{l'\Bigl(g^{-1}(|\a|^2)\Bigr)}, $$ 
and the closed set 
$$ 
K_\delta:=\Bigl\{(\a,\x) \in \Rd \times \Sd \; | \; |\a|^2 
\in g[\frac{\delta}{2},2] \Bigr\}.  
$$ 
Note that $M$ is continuous on $K_\delta$ for all $\delta>0.$ 
\hfill\break

When $\x \in \Sd$ we next denote by $\nabla^\x_{\Sd} c$ the tangential
gradient of $\x \rightarrow c(\x,\y)$ at $\x \in \Sd$ and let $T_\x$ be
the tangent space to $\Sd$ at $\x \in \Sd.$
\begin{lemma}\label{l03.08.4} Assume that $c$ is given by
(\ref{e1.5.1}) such that $l$ satisfies (\ref{e1.5.2}). Assume that $g$
is monotone on $(0,2].$ If  $0\not=\a \in \Rd,$ $\x, \y \in \Sd$ and
$\nabla^\x_{\Sd} c(\x,\y)=\a$ then $\y=M(\a,\x).$ 
\end{lemma} 
\proof{} Recall that if $\x \in \Sd$ then the orthogonal projection of
$\y \in \Rd$ onto $T_\x$ is $\y_{\|}= \y - (\x \cdot \y)\x.$ Hence if
$\y \in \Sd$ then 
\begin{equation}\label{e03.09.3} 
|\y_{\|}|^2 + (\x
\cdot \y)^2=1
 \end{equation} 
Setting $\a:=\nabla^\x_{\Sd} c(\x,\y)$ yields that 
\begin{equation}\label{e03.09.1bis} 
\a=-l'(\frac{|\x-\y|^2}{2})\y_{\|}. 
\end{equation} 
This, together
with (\ref{e03.09.3}) and the fact that 
$1- \x \cdot \y=\frac{|\x-\y|^2}{
2}$ for $\x, \y \in \Sd,$ yields that 
$$|\a|^2=\biggl(l'\Bigl(\frac{|\x-\y|^2}{2}\Bigr) 
\biggr)^2 (1- (\x \cdot\y)^2)= g(1-\x \cdot \y).  $$ Thus, 
\begin{equation}\label{e03.09.5} 
1-\x \cdot \y= g^{-1}(|\a|^2).  
\end{equation} 
We use
(\ref{e03.09.1bis}), (\ref{e03.09.5}) and the fact that 
$1- \x \cdot\y=\frac{|\x-\y|^2}{2}$ to conclude that 
$$ \y= (\x \cdot \y) \x +
\y_{\|}= \Bigl(1-g^{-1}(|\a|^2) \Bigr) \x -\frac{\a}{ 
l'\Bigl(g^{-1}(|\a|^2)\Bigr)}=M(\a,\x).$$ 
\endproof
\begin{thm}\label{t1.5} Assume that $c$ is given by (\ref{e1.5.1})
with $l \in C^2(0,+\infty)$ and  satisfying (\ref{e1.5.2}). Assume
 that $g$ given in (\ref{e03.10.1}) is monotone. Let $\mu, \nu \in
 \calP(\Sd)$ be Borel measures on the sphere $\Sd$ that vanish on
$(d-1)$--rectifiable sets. Then

(i) there exists a unique measure $\gamma_o$ that minimizes $\cost$
over $\Gam$. 

(ii) There exists a unique map $T_o: \Sd \rightarrow \Sd$ that
minimizes $\calI$ over $\calT.$ Furthermore the essential 
infimum $\inf_{\x} |T_o\x -\x|>0.$ The $\gamma_o$ is 
uniquely determined and coincides with the 
measure $(\id \times T_o)_\# \mu.$  

(iii) The map $T_o$ is invertible except on a set whose $\mu$ 
measure is null.  
\end{thm} 
\proof{} \underline{(i): Existence of the measure $\gamma_o$.} 
By proposition \ref{p03.08.2}, there exists
$\gamma \in \Gam$ and $\epsilon>0$ such that
\begin{equation}\label{e03.10.2} 
|\x-\y| \geq \epsilon 
\end{equation}
for all $(\x,\y) \in \spt \gamma.$ Since $|\x-\y| \leq 2$ on $\Sd$,
(\ref{e03.10.2}) and the fact that $c$ satisfies \smoothness  ensures
that $c$ is uniformly bounded on $\spt \gamma.$ This proves that
$\cost[ \gamma]<+\infty.$ By proposition \ref{p1.2} there is at least
one optimal measure $\gamma_o \in \Gam$ and there is a $c$--cyclically
monotone set $S \subset \Rdd$ containing the support of all optimal
measures in $\Gam.$ By remark \ref{r1.2bis} the support of $\gamma_o$
and hence $S$ contains pairs $(\x_o, \y_o),$ $ (\bar \x_o, \bar \y_o)$
such that $\bar \x_o, \bar \y_o \not \in \{\x_o , \y_o\}.$ Lemma
\ref{p1.4} ensures existence of $\delta>0$ and a $c$--concave function
$\psi$ such that $S$ is contained in the $c$--superdifferential of
$\psi$ and 
\begin{equation}\label{e.66ter09.1} 
\psi(\x)= \inf_{\y }
\Bigl\{ c(\x,\y) - \psi^c(\y) \; | \; |\x -\y| \geq \delta \Bigr\},
\quad \quad (\x \in \X).  
\end{equation} 
Clearly $\psi$ is continuous
and so is $\psi^c.$ Hence, $\partial^c \psi(\x), \partial^c \psi^c(\y)
\not = \emptyset$ for all $\x, \y \in \Sd$ and
\begin{equation}\label{e.66ter09.3} 
(\x, \y) \in \partial^c \psi \quad
\Longleftrightarrow \quad  \psi(\x)+\psi^c(\y)=c(\x,\y) 
\end{equation}
We use again (\ref{e.66ter09.1}) and the assumption that $l \in
C^2(0,\infty)$  
to conclude that $\psi$ is semiconcave on $\Sd$ and so, it is
differentiable everywhere except for a set $N \subset \Sd$ which is
$(d-1)$--rectifiable.

\hfill\break
\underline{Proof of (ii): existence and 
uniqueness of $T_o$; $\gamma_o=(\id \times T_o)_\# \mu$; 
uniqueness of $\gamma_o.$ } 
Let $\x \in \Sd \setminus N$ and $\y \in \partial^c \psi(\x)$. We use
the fact that the function $\t \rightarrow \psi(\t)+\psi^c(\y)-c(\t,\y)$  
is
differentiable and  attains its maximum at $\x$ to conclude that
\begin{equation}\label{e03.09.6} 
\nabla_{\Sd} \psi(\x)=\nabla_{\Sd}^\x
c(\x,\y).  
\end{equation} 
This, together with the Lemma \ref{l03.08.4},
implies that $\y=M(\nabla_{\Sd} \psi(\x),\x).$ Since 
\begin{equation}\label{e05.11.9}
|\x-\y| \geq \delta,
\end{equation}  
we obtain that $(\nabla_{\Sd} \psi(\x),\x) \in K_\delta.$ Hence
the map $T_o$ defined by 
\begin{equation}\label{e03.10.3}
T_o\x:=M(\nabla_{\Sd} \psi(\x),\x) 
\end{equation} is a Borel map. Note that  by (\ref{e05.11.9}), the 
map $T_o$ satisfies 
$$
|T_o\x -\x| \geq \delta 
$$ 
for all $\x \in \Sd \setminus N.$ 

Because $\mu$ vanishes on $N$ we have that $\gamma_o[N \times \Sd]=0$
and so, $T_o$ is defined $\mu$ almost everywhere. Note that we have
proved that 
\begin{equation}\label{e03.10.3bis} \partial^c \psi
\setminus (N \times \Sd) \subset {\rm graph} T_o.  
\end{equation} 

Since $\gamma_o[N \times \Sd]=0$ and  $\gamma_o \in \Gam$, we 
obtain that 
\begin{equation}\label{e03.09.7} 
\int L(\x,\y)d\gamma(\x,\y)=\int L(\x,T_o\x)d\gamma_o(\x,\y)
=\int_{\Sd}L(\x,T_o\x)d \mu(\x), 
\end{equation} 
for all $L \in
C(\Rdd).$ By (\ref{e03.09.7}) we have that 
$$ \gamma_o= (\id \times
T_o)_\# \mu, \quad T_\# \mu=\nu.  $$ 
This proves that $\gamma_o$ is
uniquely determined. We use (\ref{e03.09.7}) with $c=L$ to obtain that
$\cost[\gamma_o]=\calI[T_o].$ By (\ref{relaxation}) we conclude that
$T_o$ minimizes $\calI$ over $\calT.$

Furthermore, if $T_1$ is another minimizer of $\calI$ over $\calT$ then
$\gamma_o= (\id \times T_1)_\# \mu$ and $(\x, T_1\x) \in \spt \gamma_o
\subset \partial^c \psi$ for $\mu$ almost every $\x.$ By
(\ref{e03.10.3bis}) we have that  $T_1(\x)=T_o \x$ for $\mu$ almost
every $\x.$ This proves that $T_o$ is uniquely determined.

\hfill\break
\underline{Proof of (iii): invertibility  of the map $T_o$.}
The analogue of (\ref{e.66ter09.1}) for $\psi^c$ gives that $\psi^c$ is
semiconcave and so, the set $\tilde N$ where $\psi^c$ is not
differentiable is $(d-1)$--rectifiable. Substituting $\mu$ by $\nu$ ,
the above reasoning yields that the map $$ S_o\y= M(\nabla_{\Sd}
\psi^c(\y),\y), \quad \quad (\y \in \Sd \setminus \tilde N), $$ is such
that  $S_{o\#} \nu =\mu,$ and 
\begin{equation}\label{e03.10.3ter}
\partial^c \psi^c \setminus (\tilde N \times \Sd) \subset {\rm graph}
S_o.  
\end{equation} 
We use (\ref{e03.10.3bis}), (\ref{e03.10.3ter})
and the fact that $\partial^c \psi=\partial^c \psi^c$ to conclude that
$\id =S_o \circ T_o$ on $\Sd \setminus (N \cup S_o^{-1}(\tilde N).$
Since $\tilde N$ is $(d-1)$--rectifiable and  $S_{o\#} \nu =\mu$ we have
that $$ \mu [S_o^{-1}(\tilde N)]=\nu[\tilde N]=0.  $$ Thus, $T_o$ is
invertible on $\Sd$ up to a set whose $\mu$ measure is null.  
\endproof
\begin{rmk} If it is assumed only that $l \in C^1(0,\infty)$ and 
$\mu$ is absolutely
continuous with respect to the standard measure on $\Sd$ then
it is easy to see that $\psi, \psi^c$ are differentiable a.e.  on $\Sd$. Indeed,
it follows from  (\ref{e.66ter09.1}), (\ref{e.66ter09.3}) and
the assumption that $l \in C^1(0,\infty)$ that $\psi, \psi^c$
are locally Lipschitz. Then by Rademacher's theorem 
$\psi, \psi^c$ are differentiable a.e. Note, however, that this result
is weaker 
than the one established in Theorem \ref{t1.5}. 
\end{rmk}
%
%
%
%
%
%
%
%
%
 
%
\end{document}